\documentclass[amssymb,11pt]{amsart}
\usepackage{amscd,amssymb,euscript,amsthm}
\theoremstyle{plain}
\newtheorem{thm}{Theorem}[section] 

\newtheorem{prop}[thm]{Proposition}
\newtheorem{lem}[thm]{Lemma}
\theoremstyle{definition}
\newtheorem{defn}[thm]{Definition}
\theoremstyle{remark}
\newtheorem{rem}[thm]{Remark}

\numberwithin{equation}{section}
\newcommand{\esup}{\operatorname{ess\ sup}}
\newcommand{\comment}[1]{}
\def\<{\left<}
\def\>{\right>}
\def\cstar{$C^*$-algebra}
\subjclass[2000]{Primary 46L07; Secondary 46L52}
\begin{document}
\title{The noncommutative Choquet boundary}
\author{William Arveson}
%
%
%
%
\date{18 February, 2007}

\begin{abstract}  Let $S$ be an operator system -- a self-adjoint linear subspace 
of a unital \cstar\ $A$ 
such that $\mathbf 1\in S$ and $A=C^*(S)$ is generated by $S$.  
A {\em boundary representation} for $S$ is an 
irreducible representation $\pi$ of $C^*(S)$ on a Hilbert space 
with the property that $\pi\restriction_S$ has a unique completely 
positive extension to $C^*(S)$.  The set $\partial_S$ of all (unitary equivalence 
classes of) boundary representations is the noncommutative counterpart 
of the Choquet boundary of a function system $S\subseteq C(X)$ that separates 
points of $X$.  

It is known that the closure of the Choquet boundary of a function system $S$ 
is the \v Silov boundary of $X$ relative to $S$.  
The corresponding noncommutative problem of whether 
every operator system has  
``sufficiently many" boundary representations 
was formulated in 1969, but has 
remained unsolved despite progress on 
related issues.  In particular, it was unknown if $\partial_S\neq\emptyset$ 
for generic $S$.  
In this paper we show that every separable  
operator system has sufficiently many boundary representations.  
Our methods use separability in an essential way.  
\end{abstract}

\maketitle

\maketitle

\section{Introduction}\label{S:in}

As pointed out above, boundary representations are the 
noncommutative counterparts of points in the Choquet boundary 
of a function system in $C(X)$.  
The original motivation for introducing boundary representations 
in \cite{arvSubalgI} 
was two-fold: to provide {\em intrinsic} invariants for operator systems that could be 
calculated for specific examples, and to provide a context for showing that 
the noncommutative \v Silov boundary 
exists in general.  The first goal
was achieved in \cite{arvSubalgII} in which several concrete examples were 
worked out and applications to 
operator theory were developed - see 
Remark \ref{inRem1} for a typical example.  However, the existence of boundary representations 
and the \v Silov boundary 
was left open in general.  Subsequently, Hamana 
was able to establish the existence of the noncommutative \v Silov boundary by making 
use of his theory of injective envelopes \cite{hamInj1}, \cite{hamInj2}.  Hamana's work made 
no reference to boundary representations, and left untouched the question of their 
existence.  

More recently, Muhly and Solel \cite{muhSol1} obtained significant results 
about boundary representations in an algebraic context, and 
Dritschel and McCullough \cite{MR2132691} took a major step forward by 
showing that every unit-preserving 
completely positive map of an operator system into $\mathcal B(H)$ can be dilated 
to a completely positive map with the unique extension property 
(see Definition \ref{suDef0}).  
Significantly, that provided a new proof of the existence of the 
noncommutative \v Silov 
boundary that makes no use of injectivity.  The latter 
authors drew motivation from previous work of Agler on a model theory 
for representations of non self-adjoint operator algebras (see \cite{AgMod}
and references therein).  
On the other hand, they point out 
that their results seem to provide  
no information about the existence of boundary representations (this is 
discussed more fully in Remark \ref{inRem2}).  

An {\em operator system} is a self-adjoint linear subspace $S$ of a 
unital \cstar\ that contains the unit; we usually 
require that the \cstar\ be generated by $S$, and 
express that by writing $S\subseteq C^*(S)$. 
Let $\{\sigma_x:x\in A\}$ be a set of irreducible representations 
of $C^*(S)$.  We say that $\{\sigma_x: x\in A\}$ is 
{\em sufficient} for $S$ if 
$$
\|a\|=\sup_{x\in A}\|\sigma_x(a)\|,\qquad a\in S, 
$$
with similar formulas holding throughout the matrix hierarchy over $S$ 
in the sense that for every $n\geq 2$ and every $n\times n$ matrix 
$(a_{ij})\in M_n(S)$, we have 
\begin{equation}\label{xbEq0}
\|(a_{ij})\|=\sup_{x\in A}\|(\sigma_x(a_{ij}))\|.  
\end{equation}
If the  set of all boundary representations for $S$
is sufficient in this sense, we say that {\em $S$ has sufficiently 
many boundary representations}.  

In Theorem 2.2.3 of \cite{arvSubalgI} it was shown that, 
in all cases in which there are sufficiently many boundary representations, 
the \cstar\ $C^*(S)$ contains 
a largest closed two-sided ideal $K$ with the property that 
the natural projection $a\in C^*(S)\mapsto \dot a\in C^*(S)/K$ restricts to 
a completely isometric map on $S$;  
and in such cases one has $K=\cap\{\ker\sigma:  \sigma\in\partial S\}$.  This ideal $K$ was 
called the \v Silov boundary ideal for $S$ in \cite{arvSubalgI}.   
More recently the term has been 
contracted to {\em \v Silov ideal}, and the 
corresponding embedding $\dot S\subseteq C^*(S)/K$ 
has come to be known as the {\em $C^*$-envelope of $S$}.  As we have already 
pointed out, Hamana's work implies that the 
\v Silov ideal $K$ exists in general, independently of the 
existence of boundary representations.  
Thus, the assertion that there are sufficiently 
many boundary representations for an operator system is equivalent 
in general 
to the assertion that the \v Silov ideal is the intersection of 
the kernels of all boundary representations.  This 
is the proper noncommutative formulation 
of the statement that for every 
function system $S\subseteq C(X)$ that separates points of 
a compact Hausdorff space $X$, 
the closure of the Choquet boundary of $X$ (relative to $S$) is the \v Silov boundary 
-- the smallest closed subset of $X$ on which every function in $S$ 
achieves its norm.  

Given the central role of the 
Choquet boundary in potential theory and other parts 
of commutative analysis, it is natural to expect further applications 
of its noncommutative generalization in the future.  There 
has been a renewal of interest 
in the noncommutative \v Silov boundary, beginning around 1999 with 
work of Blecher \cite{blechSilov}, and as we have already pointed out in the preceding 
paragraphs, those developments have been fruitful.  Further 
results in these directions and additional references can be found in 
the monographs of Paulsen \cite{paulsenBk2} and Blecher and Le Merdy \cite{BlLeMbook}.  

Partly because of the promise of such developments, we were 
encouraged to return to the problem of the 
existence of boundary representations in general.  
In this paper we show that every {\em separable} operator 
system $S$ has sufficiently many boundary representations.  That is accomplished 
by first refining the theorem of Dritschel-McCullough appropriately for 
separable operator systems.  We then show that, 
given a separable Hilbert space $H$ and a 
UCP map $\phi:S\to\mathcal B(H)$ with the unique extension property, 
{\em every} direct integral decomposition of $\phi$ into 
irreducible maps gives rise to a bundle of UCP maps 
$\{\phi_x: x\in X\}$ such that $\phi_x$  is a boundary representation 
for almost every $x\in X$ with respect to the ambient measure.

The main results are Theorems \ref{dcThm1}, \ref{xbThm1} and \ref{psThm1}.  
There is further discussion of methodology and open problems 
in Section 
\ref{S:cr}.

\begin{rem}[An application of boundary representations]\label{inRem1}
The ``intrinsic" nature of the invariants associated with 
boundary representations is best illustrated by  
an example from \cite{arvSubalgII}.  If $a$ and $b$ are two irreducible compact operators with 
the property that the map 
$\lambda\mathbf 1+\mu a\mapsto \lambda\mathbf 1+\mu b$ , $\lambda, \mu\in\mathbb C$, 
is completely isometric, then $a$ and $b$ are unitarily equivalent.  Thus, 
{\em an irreducible compact operator $a$ is completely determined up to unitary 
equivalence by the internal properties of the 
two-dimensional operator space }
$$
S=\{\lambda\mathbf 1+\mu a: \lambda,\mu\in\mathbb C\}.
$$
Indeed, in \cite{arvSubalgII}, it is shown that the identity 
representation of such an 
$S$ is a boundary representation, and that any completely 
isometric map of operator systems must implement a bijection of the boundary 
representations of one operator system to those of the other.  
From these results it follows that the map 
$\lambda\mathbf 1+\mu a\mapsto 
\lambda\mathbf 1+\mu b$  
extends uniquely to a $*$-homomorphism of \cstar s.  One now 
deduces the above assertion 
from the familiar fact that an irreducible representation of the 
\cstar\ of compact operators is 
implemented by a unitary operator (Cor. 2 of Theorem 1.4.4 of \cite{arvInv}).  
\end{rem}

\begin{rem}[Terminology]\label{inRem2}
We caution the reader that in \cite{MR2132691}, the term 
{\em boundary representation} refers 
rather broadly to arbitrary UCP maps with the unique extension
property.  In this paper we adhere to the original terminology of 
\cite{arvSubalgI} and \cite{arvSubalgII}, in which {\em boundary representation} 
refers to an {\em irreducible} representation of $C^*(S)$ whose restriction to 
$S$ has the unique extension property. 
These are the objects that generalize points of the Choquet boundary 
of a function system in $C(X)$ and peak points of function algebras.  
In particular, while the results of \cite{MR2132691} 
show that there is an abundance of maps with the unique extension property, 
they provide no information about the existence of 
boundary representations in our sense of the term.  
\end{rem}

After the first version of this paper was circulated, we learned that Marius Junge has 
shown in ongoing unpublished work that every subhomogeneous 
operator system $S\subseteq \ell^\infty(M_n)$, 
$M_n$ denoting the algebra of $n\times n$ matrices, has sufficiently many 
boundary representations.  While that follows from theorem \ref{xbThm1} below when 
$S$ is separable, Junge does not assume separability.

Finally, it is with pleasure that I acknowledge valuable conversations 
 during the fall of 2002
with Narutaka Ozawa, who completely understood 
the first version of \cite{MR2132691} when I did not.  Those conversations 
led to an unpublished exposition of the results of 
\cite{MR2132691} in \cite{arvUnExt}.  In particular, Lemma \ref{suLem1} 	below
was inspired by an observation of Ozawa.  Without the paper 
\cite{MR2132691} or Ozawa's visit to Berkeley, this paper would most likely not exist.  
I would also like to thank the anonymous referee, whose careful   
reading of the manuscript led to many 
perceptive comments 
that significantly improved the readability of this paper.

\section{Maximal UCP maps of separable operator systems}\label{S:su}

In this section we discuss the unique extension property  and 
maximality for operator-valued completely positive maps 
of operator systems, and we prove a refinement of a result 
of \cite{MR2132691} that will be used below.  

We consider 
unital completely positive (UCP) maps $\phi: S\to\mathcal B(H)$, 
that is, completely positive maps that carry the unit of 
$S$ to the identity operator of $\mathcal B(H)$.  Such 
maps satisfy $\phi(x^*)=\phi(x)^*$, $x\in S$.  A 
linear map $\phi: S\to\mathcal B(H)$ that preserves the 
unit is completely positive iff it is completely contractive.  
We also recall that the more general theory of unital operator spaces (with 
unital complete contractions as maps) can be absorbed into 
the theory of operator systems (with UCP maps) because of 
the following result: 
If $S$ is a linear subspace of $C^*(S)$ containing $\mathbf 1$, 
then every completely contractive unital map of $S$ extends 
uniquely to a UCP map of $S+S^*$.  For these basic facts see 
Propositions 1.2.8--1.2.11 of  
\cite{arvSubalgI}.  

As pointed out in Remark \ref{inRem2}, 
Dritschel and McCullough have used the term {\em boundary representation} 
for UCP maps $\phi: S\to\mathcal B(H)$ 
that have unique completely positive extensions to 
representations of $C^*(S)$.  In order to avoid 
conflict in terminology, in this paper we describe 
that property as follows:  
\begin{defn}\label{suDef0}
A UCP map $\phi: S\to\mathcal B(H)$ is said to have 
the {\em unique extension property} if 
\begin{enumerate}
\item[(i)] 
$\phi$ has a unique completely positive extension 
$\tilde\phi: C^*(S)\to\mathcal B(H)$, and 
\item[(ii)] $\tilde\phi$ is a representation 
of $C^*(S)$ on $H$.  
\end{enumerate}
\end{defn}
The unique extension 
property for $\phi: S\to\mathcal B(H)$ is equivalent 
to the assertion that every extension 
of $\phi$ to a completely positive map $\phi: C^*(S)\to\mathcal B(H)$ 
should be multiplicative on $C^*(S)$.  
If the extension $\tilde \phi$ of such a map $\phi$ to $C^*(S)$ is an 
irreducible representation then the extension is a 
boundary representation in the sense of \cite{arvSubalgI}; 
otherwise it is not.

Given an operator system $S\subseteq C^*(S)$ 
and two UCP maps $\phi_k: S\to \mathcal B(H_k)$, $k=1,2$, 
we write $\phi_1\preceq \phi_2$ if $H_1\subseteq H_2$ and 
$P_{H_1}\phi_2(x)\restriction_{H_1}=\phi_1(x)$, $x\in S$;
in this event $\phi_2$ is called a {\em dilation} of $\phi_1$ and 
$\phi_1$ is called a {\em compression} of $\phi_2$.  
The relation $\preceq$ is transitive, and 
$\phi_1\preceq \phi_2\preceq \phi_1$ 
iff $H_1=H_2$ and $\phi_1=\phi_2$.  Thus, 
$\preceq$ defines a 
{\em partial ordering} of UCP maps of $S$.  Every UCP map 
$\phi: S\to\mathcal B(H)$ can be dilated in 
a trivial way by forming a direct sum 
$\phi\oplus\psi$ where $\psi: S\to\mathcal B(K)$ is 
another UCP map.  

\begin{defn}\label{def1}
A UCP map $\phi: S\to\mathcal B(H)$ 
is said to be {\em maximal}
if it has no nontrivial dilations: 
$\phi\preceq \phi^\prime \implies \phi^\prime=\phi\oplus\psi$ 
for some UCP map $\psi$.   
\end{defn}

Equivalently, $\phi$ is maximal iff for every dilation $\phi^\prime: S\to\mathcal B(H^\prime)$ 
of $\phi$ acting on $H^\prime\supseteq H$, one has $\phi^\prime(S)H\subseteq H$.  
A dilation $\phi_2$ of $\phi_1$ need not satisfy 
$H_2=[C^*(\phi_2(S))H_1]$, 
$C^*(\phi_2(S))$ denoting the \cstar\ generated 
by $\phi_2(S)\subseteq \mathcal B(H_2)$, 
but it can always be replaced 
with a smaller dilation of $\phi_1$ that
has the property, and in that case we can assert that 
$\phi_1$ is maximal iff 
the only dilation $\phi_2\succeq \phi_1$ that satisfies 
$H_2=[C^*(\phi_2(S))H_1]$ is $\phi_2=\phi_1$ itself.  
More generally, this 
reduction imposes 
an upper bound on the dimension 
of $H_2$ in terms of the dimension of $H_1$ and 
the cardinality of $S$; in particular, if $H_1$ is separable and 
$S$ is a separable operator system, then $H_2$ must be separable.

\begin{rem}[Separably acting maximal dilations]\label{suRem1}
A UCP map of an operator system 
$\phi: S\to\mathcal B(H)$ is said to be {\em separably acting} 
if $H$ is a separable Hilbert space.  Let $S$ be a separable 
operator system and let $\phi: S\to\mathcal B(H)$ be a separably 
acting UCP map.  The preceding paragraph implies that  every maximal dilation 
$\phi^\prime: S\to\mathcal B(H^\prime)$ of $\phi$ can be decomposed into a direct sum 
of maps $\phi^\prime=\tilde\phi\oplus \lambda$  where $\tilde\phi$ is 
a {\em separably acting} maximal dilation of $\phi$ and $\lambda$ 
is another UCP map.  
\end{rem}

We make repeated use of the following adaptation of 
a result of Muhly and Solel \cite{muhSol1} 
that connects maximality to the unique extension 
property:  
\begin{prop}\label{suProp1}
A UCP map $\phi: S\to\mathcal B(H)$ has the unique extension property iff it is maximal.  
\end{prop}

\begin{proof}Assume first that $\phi$ is maximal and 
let $\tilde\phi: C^*(S)\to\mathcal B(H)$ 
be a completely positive extension of it.  We have 
to show that $\tilde\phi$ is multiplicative.  By Stinespring's 
theorem, there is a representation $\sigma: C^*(S)\to\mathcal B(K)$ 
on a Hilbert space $K\supseteq H$ such that $\tilde\phi(x)=P_H\sigma(x)\restriction_H$, 
$x\in C^*(S)$.  We can 
assume that the dilation is minimal in 
that $K=[\sigma(C^*(S))H]=[C^*(\sigma(S))H]$.  By maximality of $\phi$,
$K=H$ and $\tilde\phi=\sigma $ is multiplicative.

Conversely, suppose that $\phi$ has the unique extension 
property and let $\tilde\phi: S\to \mathcal B(K)$ be a dilation of $\phi$
acting on $K\supseteq H$ with $K=[C^*(\tilde\phi(S))H]$.  
We show that $K=H$ and $\tilde\phi=\phi$.  
By Theorem 1.2.9 of \cite{arvSubalgI} $\tilde\phi$ can be extended to a completely 
positive linear map 
$\psi: C^*(S)\to \mathcal B(K)$.  Since the compression of 
$\psi$ to $H$ defines a completely positive 
map of $C^*(S)$ to $\mathcal B(H)$ that restricts 
to $\phi$ on $S$, 
the unique extension property implies that $P_{H}\psi P_{H}$ 
is multiplicative on $C^*(S)$.  So
for $x\in C^*(S)$, 
$$
P_{H}\psi(x)^*P_{H}\psi(x)P_{H}=
P_{H}\psi(x^*x)P_{H}\geq
P_{H}\psi(x)^*\psi(x)P_{H},  
$$
by the Schwarz inequality;  
hence $|(\mathbf 1-P_H)\psi(x)P_H|^2\leq 0$.  This 
implies that $H$ is invariant 
under the set of operators 
$\psi(C^*(S))\supseteq \tilde\phi(S)$, and therefore 
under $C^*(\tilde\phi(S))$.  
Both $K=[C^*(\tilde\phi(S))H]=H$ and $\tilde\phi=\phi$ follow.  
\end{proof}

We require the following refinement of the
main result of \cite{MR2132691} for separable operator systems 
and separably acting UCP maps:

\begin{thm}\label{suThm1}
Let $S\subseteq C^*(S)$ be a separable operator system, let 
$H_0$ be a separable Hilbert space,  and 
let $\phi_0:S\to \mathcal B(H_0)$ be a UCP map.  Then $\phi_0$ can 
be dilated to a separably acting UCP map $\phi:S\to\mathcal B(H)$ 
with the unique extension property.
\end{thm}

\begin{rem}[Coherent sequences of UCP maps]\label{suRem2}
In the proof of Theorem \ref{suThm1} we 
make repeated use of the following observation.  
Let $H_1\subseteq H_2\subseteq\cdots$ be an increasing sequence of 
subspaces of a  Hilbert space $H=\overline{\cup_nH_n}$.  Given a 
sequence of UCP maps $\phi_n:S\to \mathcal B(H_n)$ that is {\em coherent} 
in the sense that $\phi_n\preceq \phi_{n+1}$ for $n\geq 1$, there is a 
unique UCP map $\phi:S\to\mathcal B(H)$ that satisfies 
$\phi_n\preceq \phi$ for all $n$.  The proof is a 
straightforward exercise based on the following elementary fact: For 
every sequence of operators $a_n\in\mathcal B(H_n)$  such 
that $\sup_n\|a_n\|<\infty$ and 
$a_n=P_{H_n}a_{n+1}\restriction_{H_n}$, $n\geq 1$, 
there is a unique operator $a\in\mathcal B(H)$ satisfying 
$P_{H_n}a\restriction_{H_n}=a_n$, $n=1,2,\dots$.  Moreover, 
one obviously has $\|a\|=\sup_n\|a_n\|$. 
\end{rem}

Let $\phi: S\to \mathcal B(H)$ be a UCP map and let 
$F$ be a subset of $S\times H$.  We will 
say that $\phi$ is {\em maximal on $F$} if for every dilation 
$\psi$ of $\phi$ acting on $K\supseteq H$, we 
have 
$$
\psi(x)\xi =\phi(x)\xi,\qquad (x,\xi)\in F.  
$$
A UCP map $\phi:S\to\mathcal B(H)$ is maximal iff 
it is maximal on $S\times H$.  
Note that if $\phi$ is maximal on 
$F\subseteq S\times H$ and $\psi\succeq\phi$, 
then $\psi$ is maximal on $F$ as well.  
Our proof of Theorem \ref{suThm1} uses the following result,  
inspired by an observation of N. Ozawa.

\begin{lem}\label{suLem1}Let $S$ be a separable operator 
system.  For every 
UCP map $\phi:S\to\mathcal B(H)$ where $H$ is a separable 
Hilbert space and 
every $(x,\xi)\in S\times H$, there is a separably-acting dilation 
of $\phi$ that is maximal on $(x,\xi)$.  
\end{lem}

\begin{proof}
Since for every dilation $\psi\succeq\phi$ we have 
$\|\psi(x)\xi\|\leq \|x\|\cdot\|\xi\|<\infty$, we can find 
a separably acting dilation $\phi_1$ of $\phi$ for which $\|\phi_1(x)\xi\|$ is as close 
to $\sup\{\|\psi(x)\xi\|: \psi\succeq\phi\}$ as we wish.  
Continuing inductively, we find a sequence of separably acting UCP 
maps $\phi\preceq \phi_1\preceq \phi_2\preceq\cdots$ such that 
$\phi_n: S\to\mathcal B(H_n)$, 
$H\subseteq H_1\subseteq H_2\subseteq\cdots$, and 
$$
\|\phi_{n+1}(x)\xi\|\geq \sup_{\psi\succeq\phi_n}\|\psi(x)\xi\| - 1/n.    
$$
Let $H_\infty$ be the closure of the union $\cup_n H_n$ and let 
$\phi_\infty: S\to \mathcal B(H_\infty)$ be the unique UCP map 
that compresses to $\phi_n$ on $H_n$ for every $n$ (see Remark \ref{suRem2}).  Note that 
$\phi_\infty$ is maximal on $(x,\xi)$.  Indeed, if $\psi\succeq\phi_\infty$ 
then $\psi\succeq\phi_{m}$ for every $m\geq 1$. Fixing $n\geq 1$ and  $\xi\in H_{n+1}$, 
choose $m>n$.  Then   
$$
\|\phi_\infty(x)\xi\|\geq\|P_{H_{m+1}}\phi_\infty(x)\xi\|=
\|\phi_{m+1}(x)\xi\|\geq \|\psi(x)\xi\|-1/m,   
$$
so that $\|\phi_\infty(x)\xi\|\geq \|\psi(x)\xi\|$ because $m$ 
can be arbitrarily large.  
Hence  
$$
\|\psi(x)\xi-\phi_\infty(x)\xi\|^2=
\|\psi(x)\xi-P_{H_\infty}\psi(x)\xi\|^2=
\|\psi(x)\xi\|^2-\|\phi_\infty(x)\xi\|^2\leq 0
$$ 
so that $\psi(x)\xi=\phi_\infty(x)\xi$.  The assertion follows 
since $\overline{\cup_nH_{n+1}}=H_\infty$.   
\end{proof}

\begin{proof}[Proof of Theorem \ref{suThm1}]
We claim first that $\phi_0$ can be dilated to a separably acting UCP map 
$\phi_1: S\to \mathcal B(H_1)$ that is maximal on $S\times H_0$.  
To that end, let $C$ be a countable dense subset of $S$,  
let $D$ be a countable norm-dense subset of $H_0$, and 
enumerate the elements of $C\times D=\{z_1,z_2,\dots\}$.  
We claim 
that there is a sequence of separably acting UCP maps 
$\omega_n: S\to\mathcal B(K_n)$, $n\geq 1$, such that 
\begin{enumerate}
\item[(i)]$\phi_0\preceq \omega_1\preceq \omega_2\preceq\cdots$, and 
\item[(ii)]$\omega_n$ is maximal on $\{ z_1,\dots,z_n\}$.  
\end{enumerate}
Indeed, Lemma \ref{suLem1} implies the existence of a 
separably acting dilation $\omega_1$ of $\phi_0$ that is maximal 
on $z_1$.  Given that $\omega_1,\dots,\omega_n$ have been defined 
and satisfy (i) and (ii), the same reasoning gives a separably 
acting dilation $\omega_{n+1}$ of $\omega_n$ that is maximal 
on $z_{n+1}$; since $\omega_{n+1}$ dilates each of the preceding 
maps, it must also be maximal on $z_1,\dots,z_n$.  
Once one is given such a sequence $\omega_1, \omega_2,\dots$, one can let 
$H_1$ be the closure of $\cup_nK_n$ and let $\phi_1$ be the unique 
UCP map of $S$ into $\mathcal B(H_1)$ that compresses to 
$\omega_n$ on each $K_n$.

By an obvious induction on the preceding fact, 
one obtains an increasing sequence of separable Hilbert spaces 
$H_0\subseteq H_1\subseteq H_2\subseteq\cdots$ and UCP maps 
$\phi_n: S\to\mathcal B(H_n)$ such that $\phi_{n+1}$ is a dilation 
of $\phi_n$ that is maximal 
on $S\times H_n$, $n=0,1,2,\dots$.  Let $H_\infty$ be 
the closure of $\cup_nH_n$ and let $\phi_\infty: S\to\mathcal B(H_\infty)$ 
the unique UCP map that compresses to $\phi_n$ on $H_n$ for every $n\geq1$.  
Note that for 
every dilation $\psi: S\to\mathcal B(K)$ of $\phi_\infty$ and every 
$n\geq 1$, both $\psi$ and $\phi_\infty$ are dilations of $\phi_{n+1}$, 
so by maximality of $\phi_{n+1}$ on $S\times H_n$ we have 
$$
\psi(x)\xi= \phi_{n+1}(x)\xi =\phi_\infty(x)\xi, 
\qquad (x,\xi)\in S\times H_n.  
$$
It follows that $\phi_\infty$ is maximal on $S\times\cup_n H_n$, 
hence on its closure $S\times H_\infty$.   
\end{proof}

\begin{rem}[Significance of the relation $\phi\preceq\psi$]
There is a weaker and considerably more flexible ordering $\prec$ of UCP maps, 
in which for UCP maps $\phi_k: S\to\mathcal B(H_k)$, $k=1,2$,   $\phi_1\prec\phi_2$ 
means that there is an isometry $V: H_1\to H_2$ such that 
$\phi_1(a)=V^*\phi_2(a)V$, $a\in S$.  Equivalently,  $\phi_1$ is {\em unitarily 
equivalent} to a map $\phi_1^\prime$ satisfying $\phi_1^\prime\preceq\phi_2$.  
While there is 
a variation of Theorem \ref{suThm1} that makes a similar assertion about the 
maximality properties of the 
relation  $\prec$, one cannot prove it by a {\em verbatim} replacement 
of $\preceq$ with $\prec$ in the 
above arguments.  This subtle 
difficulty becomes apparent when one attempts to 
associate a single ``limit" 
UCP map with a weakly increasing sequence $\phi_1\prec\phi_2\prec\cdots$ 
as we did in Remark \ref{suRem2}.  
\end{rem}

\section{Borel cross sections}\label{S:mt}

Given standard Borel spaces $X$, $Y$ and a surjective Borel map $f: X\to Y$, 
 a {\em cross section} for $f$ is a Borel function $g: Y\to X$ such 
that $f\circ g$ is the identity map of $Y$.  Not every surjective Borel 
map of standard Borel spaces has a Borel cross section.  Indeed, 
there exist Borel subsets $X$ of the unit square 
$[0,1]\times[0,1]$ with 
the property that the projection $p(x,y)=x$ maps $X$ onto $[0,1]$, 
but $p$ does not have a Borel cross section \cite{LusSel},  \cite{NovSel}; 
a simpler example was given by Blackwell in \cite{BlGraph}.  
  The following 
result provides the measure-theoretic substitute that we require.

\begin{thm}\label{mtThm1}
Let $X$, $Y$ be standard Borel spaces, let $f: X\to Y$ be a surjective Borel map, 
and let $\mu$ be a finite positive 
measure on $Y$.  Then 
there is a Borel set $N\subseteq Y$ of measure zero and a 
Borel map $g: Y\setminus N\to X$ such that $f\circ g$ is the identity 
map on $Y\setminus N$.  
\end{thm}

Theorem \ref{mtThm1} is part of the lore of the subject; but since we lack 
a convenient reference and the result is needed below, 
we briefly indicate how one deduces it 
from a selection theorem proved in \cite{arvInv}.  
Recall that a subset $A\subseteq X$ of a standard Borel space $X$ is said to be 
{\em absolutely measurable} if, for every finite positive measure 
$\mu$ on $X$, there are Borel sets $E_\mu$, $F_\mu$ such that 
$E_\mu\subseteq A\subseteq F_\mu$ and $\mu(F_\mu\setminus E_\mu)=0$.  
The class of all absolutely measurable subsets of $X$ is a 
$\sigma$-algebra containing the Borel sets.  Analytic sets are 
examples of absolutely measurable sets that need not be Borel sets 
(see Section 3.4 of \cite{arvInv}).  
A function $f:X\to Y$ from a standard Borel space $X$ to a Borel 
space $Y$ is said to be  absolutely measurable if $f^{-1}(E)$ is 
absolutely measurable for every Borel set $E\subseteq Y$.  
Theorem 3.4.3 of \cite{arvInv} specializes to the 
following assertion in this context:

\begin{thm}\label{mtThm2}
Let $X$ be a standard Borel space and let 
$Y$ be a countably separated Borel space.  Then every surjective 
Borel map $f:X\to Y$ has an absolutely 
measurable cross section.  
\end{thm}

To deduce Theorem \ref{mtThm1} from Theorem \ref{mtThm2}, note that 
for every finite measure $\mu$ on $Y$, 
an absolutely measurable cross section $g: Y\to X$ for 
$f$ must agree almost everywhere 
$(d\mu)$ with a Borel function $g_\mu: Y\to X$, where of course 
$g_\mu$ depends on $\mu$.  Indeed, this is obvious if $X$ is finite or 
countable; if $X$ is uncountable then it is Borel isomorphic to the 
unit interval $[0,1]$,  $g$ becomes a real valued function 
in $L^\infty(Y,\mu)$, and such a function must agree with a Borel 
function $g_\mu$ almost everywhere $(d\mu)$.  
Letting $N\subseteq Y$ 
be a Borel set of $\mu$-measure zero such that $g=g_\mu$ on $Y\setminus N$, 
one obtains $f\circ g_\mu(y)=f\circ g(y)=y$ for all $y\in Y\setminus N$, and Theorem 
\ref{mtThm1} follows.  

We also require the following result, which follows from Theorem 3.3.4 
of \cite{arvInv}.  Recall that a subset $A$ of a standard Borel space $X$ is said to 
be {\em analytic} if there is an analytic set in a Polish space that is 
isomorphic to $A$ with its relative Borel structure.  

\begin{thm}\label{mtThm3}
Let $X$, $Y$ be standard Borel spaces and let $f: X\to Y$ be a 
Borel map.  Then $f(X)$ is an analytic set in $Y$ and  is therefore  
absolutely measurable.  
\end{thm}

\section{Borel Families of UCP maps}\label{S:fa}
In this section we discuss families of UCP maps of $S$ and 
their basic measurability properties.  
Throughout, $S$ will denote a separable operator system and $X$ will denote 
a standard Borel space.  

Operator algebraists traditionally refer to Chapitre II of \cite{dixVn} for the 
basic disintegration 
theory of representations of \cstar s.  We prefer a different 
- though roughly equivalent - formulation that is 
based on standard Borel spaces, and which 
runs more parallel to the modern theory of vector bundles 
over topological spaces.  Thus 
we shall formulate the basic structures with some care, 
so that key results in the following sections concerning 
decomposable maps of operator spaces 
can be readily deduced from the development of \cite{arvInv}, Chapters 3 and 4.   
For clarity, we have included more generality (and more detail) 
than is actually required below.  

Let $H=\{H_x: x\in X\}$ be a standard bundle of 
separable Hilbert spaces over $X$.  
More precisely, we are given a standard Borel 
space $H$ and a surjective Borel map $p: H\to X$ with the property that 
$H_x=p^{-1}(x)$ is a separable Hilbert space for every $x$, such that vector addition, 
multiplication by scalars, 
and the inner product are Borel-measurable.   
Notice that 
both vector addition and the inner product are defined on the 
following Borel subset of $H\times H$
$$
B=\{(\xi,\eta)\in H\times H: p(\xi)=p(\eta)\}\subseteq H\times H,   
$$ 
and the requirements are that both functions 
$(\xi,\eta)\in B\mapsto \xi+\eta\in H$ and $(\xi,\eta)\in B\mapsto \langle\xi,\eta\rangle\in \mathbb C$ 
should be Borel-measurable.  Measurability of scalar multiplication 
means that the function 
$$
(\lambda,\xi)\in \mathbb C\times H\mapsto \lambda\cdot\xi\in H
$$
should be Borel.  If one is given a complex-valued Borel function $\mu:H\to\mathbb C$, 
then the function $\xi\in H\mapsto \mu(\xi)\xi\in H$ can be expressed  
as the composition of the Borel function $\xi\in E\mapsto (\mu(\xi),\xi)\in\mathbb C\times E$, 
with the Borel function $(\lambda,\xi)\in \mathbb C\times E\mapsto \lambda\xi\in E$, 
and is therefore a Borel function.  It follows  
from this observation 
that for every complex-valued Borel function $\lambda:X\to\mathbb C$
and every Borel section $x\in X\mapsto \xi(x)\in H_x$ of $p:H\to X$, the function 
$$
x\in X\mapsto \lambda(x)\xi(x)=\lambda(p(\xi(x)))\xi(x)=\lambda\circ p(\xi(x))\xi(x) 
$$ 
is another Borel section; indeed, it is the composition of the Borel function $x\mapsto \xi(x)$ 
with the above map $\xi\mapsto\mu(\xi)\xi$ in which $\mu=\lambda\circ p$.      

Though it is not 
necessary for the development, it will be 
convenient to assume that $H_x\neq\{0\}$ for every 
$x\in X$.  

The last and most important axiom for a standard Hilbert bundle is 
that there should exist  
a sequence of Borel sections $\xi_n: x\in X\mapsto \xi_n(x)\in H_x$, 
$n=1,2,\dots$, such that $H_x$ is the closed linear span of $\{\xi_1(x),\xi_2(x),\dots\}$ 
for every $x\in X$.  
There is a natural notion of (unitary) isomorphism of Hilbert bundles $p^k: H^k\to X$ over $X$, 
$k=1,2$, namely 
an isomorphism of Borel spaces $U: H^1\to H^2$ with the property that
$U$ restricts to a unitary operator from  $H^1_x$ to $H^2_x$ for every $x\in X$.

For each $n=\infty, 1,2,\dots$, let $X_n=\{x\in X: \dim H_x=n\}$.  
A straightforward - though somewhat tedious - argument shows that each $X_n$ is a Borel set.
Alternately, one can prove that by making appropriate 
use of the Gram-Schmidt procedure that we will describe momentarily.  
More significantly, the restriction of the bundle to $X_n$ is isomorphic 
to the trivial bundle $p:X_n\times \mathbb C^n\to X_n$ if $n<\infty$, 
or $p:X_\infty\times \ell^2\to X_\infty$ if $n=\infty$.  To sketch the proof of 
the latter fact, we assume 
first that $n<\infty$.  The preceding axioms,  along with the Gram-Schmidt 
procedure applied to the given sequence of sections $\{\xi_1(x), \xi_2(x),\dots\}$ 
restricted to $X_n$, 
allow one to construct 
a new sequence of Borel sections 
$$
e_i: x\in X_n\mapsto e_i(x)\in H_x,\qquad 1\leq i\leq n, 
$$
with the property that $\{e_1(x),e_2(x),\dots,e_n(x)\}$ is an orthonormal basis 
for $H_x$ for each $x\in X_n$.  In the case $n=\infty$, 
the above procedure generates a 
countably infinite Borel-measurable orthonormal basis $\{e_1(x),e_2(x),\dots\}$ for $H_x$ 
for every $x\in X_\infty$.  Once one has such an orthonormal basis for each 
$x\in X_n$ and each $n$, it becomes obvious how  
one can write down a unitary isomorphism of Hilbert bundles that 
trivializes each restricted bundle $H\restriction_{X_n}$.  For example, for 
finite $n$  
one can define a unitary operator $U_x: H_x\to \mathbb C^n$ for every $x\in X_n$ by 
$$
U_x(\xi)=(\langle\xi,e_1(x)\rangle,\cdots,\langle \xi,e_n(x)\rangle)\in\mathbb C^n,
\qquad \xi\in H_x,   
$$
and one obtains a trivializing isomorphism 
$U: H\restriction_{X_n}\to X_n\times\mathbb C^n$ as the total map $U(x,\xi)=(x,U_x\xi)$  
associated with this family of unitary operators.   

As an illustration, using the above remarks it is not hard to show that a section 
$x\in X\mapsto \eta(x)\in H_x$ is Borel-measurable iff each of the 
complex-valued inner products $x\in X\mapsto\langle\eta(x),\xi_n(x)\rangle$, 
$n=1, 2,\dots$,  is a Borel function (see Remark \ref{faRem1} for more detail).   
More generally, these remarks show that 
one can reduce the analysis of standard Hilbert bundles to the case of flat Hilbert 
bundles, and this procedure of 
``flattening and piecing together" is useful for clarifying the nature of 
more complex structures 
associated with a Hilbert bundle, as we will see momentarily.    

There is also a natural notion of bounded measurable family of operators $A=\{A_x: x\in X\}$ 
on a Hilbert bundle, namely a Borel map $A:H\to H$ that restricts to a bounded linear 
operator $A_x\in\mathcal B(H_x)$ on each fiber such that $\sup_{x\in X}\|A_x\|<\infty$.  
The set of all such operator families is a unital \cstar.  
More generally, there is a natural notion 
of a weak measurability for a family $\phi=\{\phi_x: x\in X\}$ 
of UCP maps of $S$: 

\begin{defn}\label{faDef1}
By a {\em family} of UCP maps on $S$ we mean 
set $\{\phi_x: x\in X\}$ of UCP maps $\phi_x: S\to\mathcal B(H_x)$ 
indexed by $X$ 
with the property that for every 
pair of Borel sections $\xi,\eta$ of $H$ and for every 
$s\in S$, the inner product 
$
x\mapsto\langle\phi_x(s)\xi(x),\eta(x)\rangle
$ 
defines a complex-valued Borel function on $X$.  
\end{defn}

We will write a family of UCP maps in the more descriptive way as 
$$
\phi=\{\phi_x: S\to \mathcal B(H_x): x\in X\}.
$$
A family of UCP maps is a Borel cross section  
of a bundle of UCP maps that is described as follows.  
For each $x\in X$, the set  $UCP(S,\mathcal B(H_x))$ 
of all UCP maps from $S$ into $\mathcal B(H_x)$ is a 
convex set of linear maps defined on $S$.    
The total space of this collection of maps is 
$$
UCP(X,S,\mathcal B(H))=\{(x,\phi): x\in X,\quad \phi\in UCP(S,\mathcal B(H_x)\},   
$$
with projection $p: UCP(X,S,\mathcal B(H))\to X$ given by 
$p(x,\phi)=x$, and one has $UCP(S,\mathcal B(H_x))=p^{-1}(x)$.  

There is a natural way to make $UCP(X,S,\mathcal B(H))$ 
into a standard Borel space in such a way that $p$ is a Borel map.  To see that, 
assume first that the fibers $H_x$ do not depend on $x$, say 
$H_x=K$, $x\in X$, $K$ being a separable Hilbert space.  Then for each $x\in X$ we have 
$$
UCP(S,\mathcal B(H_x))=UCP(S,\mathcal B(K)) 
$$
and therefore $p:UCP(X,S,\mathcal B(H)))\to X$ becomes the trivial 
family of sets 
$$
p:X\times UCP(S,\mathcal B(K)))\to X.  
$$
The space of maps $UCP(S,\mathcal B(K))$ 
carries a BW-topology  (see \cite{arvSubalgI}), and since 
$S$ and $K$ are both separable, this topological space is metrizable and compact.  
Hence the cartesian product $X\times UCP(S,\mathcal B(K)))$ 
is a standard Borel space and $p$ becomes a Borel map for which 
$p^{-1}(x)=UCP(S,\mathcal B(K))$.  
So in this case, $UCP(X,S,\mathcal B(H))$ is a trivial 
standard Borel bundle 
of compact convex metrizable spaces of linear maps
$p:X\times UCP(S,K))\to X$.  

If $p:H\to X$ is merely isomorphic to a trivial bundle $p:X\times K\to X$ 
of separable Hilbert spaces, then one can make use of that isomorphism of Hilbert bundles 
in the obvious way 
to transfer the Borel structure of the flat bundle $p:X\times UCP(S,\mathcal B(K)))\to X$ to 
$p:UCP(X,S,\mathcal B(H))\to X$.  Finally, in the 
general case one decomposes $X=X_\infty\cup X_1\cup X_2\cup\cdots$ 
into a disjoint union, thereby decomposing the space 
of maps  $UCP(X,S,\mathcal B(H))$ into a disjoint union, and 
one introduces a (standard) Borel structure on this disjoint 
union of standard Borel spaces in the usual way.  
Thus, this ``flattening and piecing together" procedure 
allows one to introduce a standard Borel structure on the bundle of maps 
$p:UCP(X,S,\mathcal B(H))\to X$.

\begin{rem}[Measurability of sections]\label{faRem1} 
Let $\{H_x: x\in X\}$ be a standard Hilbert bundle and 
suppose that, for each $x\in X$, we have a UCP map  
 $\phi_x: S\to \mathcal B(H_x)$. 
The function $x\in X\mapsto \phi_x\in UCP(S,\mathcal B(H_x))$ 
is a section of the bundle of maps 
$$
p:UCP(X,S,\mathcal B(H))\to X
$$ 
and Definition \ref{faDef1} makes an assertion about the measurability of 
this section that is not identical {\em verbatim} with the definition 
of measurability that accompanies the Borel structure defined on 
$UCP(X,S,\mathcal B(H))$.    
Since this is a central issue, we sketch a proof that 
the two definitions of measurability are in fact equivalent.  
Indeed, according to the 
``flattening and piecing together" procedure, it is enough to check the 
equivalence of the two definitions in the case where $\{H_x: x\in X\}$ is 
the constant Hilbert bundle $H_x=K$, $x\in X$, $K$ being a separable Hilbert space.  
A few moments' thought shows that the equivalence of the two definitions 
is a consequence of the equivalence of the following assertions about 
self-adjoint operator functions $x\in X\mapsto A_x=A_x^*\in \mathcal B(K)$: 
\begin{enumerate}
\item[(i)] For each $\xi\in K$, the function $x\in X\mapsto A_x\xi\in K$ is Borel measurable.  
\item[(ii)]For each $\xi,\eta\in K$, the complex-valued function 
$x\in X\mapsto\langle A_x\xi,\eta\rangle$ 
is Borel measurable.  
\item[(iii)] For every pair of Borel functions $\xi, \eta: X \to K$, 
the complex-valued function $x\mapsto \langle A_x\xi(x),\eta(x)\rangle$ is Borel measurable.  
\end{enumerate}
\end{rem}
Indeed, the equivalence of (i) and (ii) follows from the fact that $K$ is a Polish 
space relative to its norm topology, and that the sigma algebra generated by 
the weak topology of $K$ coincides with the sigma algebra generated by its norm topology 
(for example, see Theorem 3.3.5 of \cite{arvInv}).  
Obviously (iii)$\implies$(ii), and (ii)$\implies$(iii) 
becomes obvious as well after one expands the inner product of 
(iii) using an orthonormal basis $e_1, e_2,\dots$ for $K$, 
\begin{align*}
\langle A_x\xi(x),\eta(x)\rangle&=\sum_n\langle A_x\xi(x),e_n\rangle\langle e_n,\eta(x)\rangle
=\sum_n\langle \xi(x),A_xe_n\rangle\langle e_n,\eta(x)\rangle
\\
&=
\sum_{m,n}\langle\xi(x),e_m\rangle\langle e_m,A_xe_n\rangle\langle e_n,\eta(x)\rangle.  
\end{align*}

\section{Sections of the dilation bundle}\label{S:db}

Throughout this section,  $\{H_x:x\in X\}$ will denote a standard Borel  
Hilbert bundle and $\phi=\{\phi_x: S\to\mathcal B(H_x): x\in X\}$ 
will denote a family of UCP maps defined on a separable operator system $S$.  
We introduce the bundle of dilations of the family $\phi$.  Given a probability 
measure $\mu$ on $X$, we show 
that nontrivial Borel-measurable sections of the dilation bundle exist 
whenever $\phi_x$ fails to be maximal for $\mu$-almost every $x\in X$.

\begin{defn}\label{dbDef1}
By a {\em dilation} of the family $\phi$
we mean a Borel family of UCP maps 
$\psi=\{\psi_x: S\to\mathcal B(H_x\oplus\ell^2): x\in X\}$ that compresses 
pointwise to $\phi$ in the sense that 
$$
\langle\psi_x(s)\xi,\eta\rangle=\langle\phi_x(s)\xi,\eta\rangle, 
\qquad s\in S, \quad \xi, \eta\in H_x,\quad x\in X.  
$$
\end{defn}

Notice that the dilations of $\phi$ act on a Hilbert bundle $\{K_x: x\in X\}$ that is 
related to the Hilbert bundle $\{H_x: x\in X\}$ of $\phi$ in a particularly 
concrete way, namely  $K_x=H_x\oplus \ell^2$ for 
all $x\in X$.  

The dilations of $\phi$ are the Borel sections of a 
bundle $\mathcal D=\mathcal D^\phi$ that we now define as follows.  
For each $x\in X$,  
let $\mathcal D_x$ be the set of all UCP maps $\psi: S\to\mathcal B(H_x\oplus\ell^2)$ 
that compress to $\phi$ on $H_x$: 
$$
\langle\psi(s)\xi,\eta\rangle=\langle\phi_x(s)\xi,\eta\rangle, \qquad s\in S,\quad \xi,\eta\in H_x.  
$$
The dilation bundle $p:\mathcal D\to X$ is the total space 
$$
\mathcal D=\{(x,\psi): x\in X,\quad \psi\in \mathcal D_x\}
$$  
with natural projection $p(x,\psi)=x$.  We can view $\mathcal D$ as 
a subset of the bundle of maps UCP$(X,S,\mathcal B(H\oplus\ell^2))$, and as such it 
inherits a relative 
Borel structure making $p:\mathcal D\to X$ a Borel map.  

\begin{prop}\label{dbProp0}
$\mathcal D$ is a Borel subset of $UCP(X,S,\mathcal B(H\oplus\ell^2))$, hence 
it is a standard Borel space.  
The dilations of $\phi$ are the Borel sections of the bundle 
$p:\mathcal D\to X$.  
\end{prop}

\begin{proof}
For each $x\in X$, let $P_x$ be the projection of $H_x\oplus \ell^2$ onto 
$H_x$.  For each $x\in X$ there is a natural compression map $\gamma_x$ 
of $UCP(S,\mathcal B(H_x\oplus\ell^2))$ 
onto $UCP(S,\mathcal B(H_x))$ defined by 
$$
\gamma_x(\psi)(s)=P_{H_x}\psi(s)\restriction_{H_x}, \qquad 
s\in S, 
$$
and the total map $\gamma$ is defined at the level of bundles by  
$$
\gamma: (x,\psi)\in UCP(X,S,\mathcal B(H\oplus\ell^2))\mapsto (x,\gamma_x(\psi))\in 
UCP(X,S,\mathcal B(H)).  
$$
The total map $\gamma$ is a Borel map because both 
the inclusion $H_x\subseteq H_x\oplus\ell^2$ 
and its adjoint $P_{H_x}: H_x\oplus\ell^2\to H_x$ define Borel maps of (standard) Hilbert 
bundles whose fibers are, respectively, isometries and co-isometries.  
We can now exhibit $\mathcal D$ as a Borel set 
in $UCP(X,S,\mathcal B(H\oplus\ell^2))$ by way of 
$$
\mathcal D=\{(x,\psi)\in UCP(X,S,\mathcal B(H\oplus\ell^2)): \gamma(x,\psi)=(x,\phi)\},   
$$
after noting that the inverse image of a singleton 
under a Borel map of standard Borel spaces $\gamma:U\to V$ is a Borel set in $U$.  

The second sentence now follows from Remark \ref{faRem1},  
which characterizes Borel sections of bundles of UCP maps.  
\end{proof}

In particular, note that 
for fixed $x\in X$, the fiber $\mathcal D_x=p^{-1}(x)$ is a closed convex subset of the BW-compact 
space of maps UCP$(S,\mathcal B(H_x\oplus\ell^2))$.

\begin{prop}\label{dbProp1}Let $s\in S$, let $\xi$ be a Borel section of 
$\{H_x: x\in X\}$, and fix $x\in X$.  Then $\phi_x$ is maximal on $(s,\xi(x))$ iff 
\begin{equation}\label{dbEq1}
\|\psi(s)\xi(x)\|\leq \|\phi_x(s)\xi(x)\|,\qquad \forall \psi\in \mathcal D_x.   
\end{equation}
\end{prop}

\begin{proof}  If $\phi_x$ is maximal on $(s,\xi(x))$ then 
for every dilation $\psi\in \mathcal D_x$ we have $\psi(s)\xi(x)=\phi_x(s)\xi(x)$, 
and (\ref{dbEq1}) follows with equality.  Conversely, assume that every 
dilation $\psi\in\mathcal D_x$ 
satisfies (\ref{dbEq1}), and let $\omega: S\to \mathcal B(K)$ be an 
arbitrary dilation of $\phi_x$, with $K\supseteq H_x$.  By the remarks 
following Definition \ref{def1}, $\omega$ decomposes into a direct sum 
$\omega_0\oplus\lambda$ where $\omega_0$ acts on a separable Hilbert 
space $K_0$ with $K\supseteq K_0\supseteq H_x$.  We can identity $K_0$ with 
either $H_x\oplus\mathbb C^k$ 
or $H_x\oplus\ell^2$ and replace $\omega_0$ with a unitarily equivalent 
map $\omega_0^\prime$ from $S$ to the respective space of operators 
such that $\phi_x\preceq\omega_0^\prime$.  
If $K_0=H_x\oplus\ell^2$ then we set $\psi=\omega_0^\prime$; 
otherwise we set $\psi=\omega_0^\prime\oplus\mu: S\to\mathcal B(H_x\oplus\ell^2)$ 
where $\mu$ is an arbitrary UCP map from $S$ into 
$\mathcal B(\ell^2\ominus\mathbb C^k)$ that fills out the difference.  Thus in all cases 
we have exhibited a map $\psi\in\mathcal D_x$ such that 
$\|\omega(s)\xi(x)\|\leq \|\psi(s)\xi(x)\|$.  From (\ref{dbEq1}) it follows that  
\begin{align*}
\|\omega(s)\xi(x)-\phi_x(s)\xi(x)\|^2&=\|\omega(s)\xi(x)-P_{H_x}\omega(s)\xi(x)\|^2
\\
&=
\|\omega(s)\xi(x)\|^2-\|P_{H_x}\omega(s)\xi(x)\|^2
\\
&\leq\|\psi(s)\xi(x)\|^2-\|\phi_x(s)\xi(x)\|^2\leq 0.    
\end{align*}
Hence $\omega(s)\xi(x)=\phi_x(s)\xi(x)$, so that $\phi_x$ is maximal on $(s,\xi(x))$.  
\end{proof}

\begin{defn}\label{dbDef2}  Fix a probability measure 
$\mu$ on $X$.  
A family of UCP maps 
$\phi=\{\phi_x:S\to\mathcal B(H_x)\}$ 
is said to be {\em maximal $\mu$-almost everywhere} if there is 
a Borel set $N\subseteq X$ such that $\mu(N)=0$ and 
$\phi_x$ is maximal for all $x\in X\setminus N$.  
\end{defn}

\begin{lem}\label{dbLem1}
Let $\phi=\{\phi_x: S\to\mathcal B(H_x): x\in X\}$ be a family 
of UCP maps of $S$ and let $\mu$ be a probability measure on $X$.
If for every $s\in S$ and every Borel section $\xi: x\in X\to H_x$ of 
$\{H_x: x\in X\}$ there is a Borel set $N_{s,\xi}\subseteq X$ such that 
$\mu(N_{s,\xi})=0$ and 
$\phi_x$ satisfies (\ref{dbEq1}) for every $x\in X\setminus N_{s,\xi}$, then 
$\phi$ is maximal $\mu$-almost everywhere.    
\end{lem}

\begin{proof}
Let $C$ be a countable subset of $S$ that is dense in 
$S$ and let $D=\{\xi_1,\xi_2,\dots\}$ be a sequence of Borel sections 
of the Hilbert bundle $\{H_x: x\in X\}$ such that $\{\xi_1(x),\xi_2(x),\dots\}$ 
is dense in 
$H_x$ for every $x\in X$.  By hypothesis, for every $(s,\xi)\in C\times D$, there is a Borel 
set $N_{s,\xi}$ of measure zero such that 
$$
\sup\{\|\psi(s)\xi(x)\|: \psi\in\mathcal D_x\}=\|\phi_x(s)\xi(x)\|,\qquad x\in X\setminus N_{s,\xi}.  
$$
Let $N$ be the (countable) union of the sets $N_{s,\xi}$,  $(s,\xi)\in C\times D$.   
Then $\mu(N)=0$, and  for every $x\in X\setminus N$ and every $\psi\in \mathcal D_x$, 
we have 
$$
\|\psi(s)\xi(x)\|\leq \|\phi_x(s)\xi(x)\|, \qquad s\in C,\quad \xi\in D.     
$$
It follows that when $x\in X\setminus N$ we have 
$\|\psi(s)\eta\|\leq \|\phi_x(s)\eta\|$ for every $s\in S=\overline C$, 
every  $\eta\in H_x=\overline{\{\xi_1(x),\xi_2(x), \dots\}}$, 
and every $\psi\in \mathcal D_x$.   Proposition \ref{dbProp1} implies that 
$\phi_x$ is maximal for every 
$x\in X\setminus N$.  
\end{proof}

 The following result provides the key step in the proof of Theorem \ref{dcThm1}:   

\begin{thm}\label{dbThm1}  Let $\mu$ be a probability measure 
on $X$ and let 
$$
\phi=\{\phi_x:S\to\mathcal B(H_x): x\in X\}
$$ 
be a family of UCP maps that fails to be maximal $\mu$-almost everywhere.  
Then there is an operator $s\in S$, a Borel 
section $\xi$ of 
$\{H_x:x\in X\}$, a Borel subset $X_0\subseteq X$ of positive 
measure and a Borel function $x\in X_0\mapsto \psi_x\in \mathcal D_x$ such that 
\begin{equation}\label{dbEq2}
\|\psi_x(s)\xi(x)\|>\|\phi_x(s)\xi(x)\|,\qquad x\in X_0.  
\end{equation}
In particular, for each $x\in X_0$, $\psi_x$ is a dilation 
of $\phi_x$ that does not decompose into a direct sum 
$\phi_x\oplus\lambda_x$.  
\end{thm}

\begin{proof}  By Lemma \ref{dbLem1}, there is an operator 
$s\in S$ and a Borel section $\xi: x\in X\mapsto \xi(x)\in H_x$ with 
the property that $\phi_x$ fails to satisfy (\ref{dbEq1}) $\mu$-almost 
everywhere.  This simply means that if 
$E\subseteq X$ is a Borel set with the property that 
$$
 \sup_{\psi\in\mathcal D_x}\|\psi(s)\xi(x)\|\leq \|\phi_x(s)\xi(x)\|
$$ 
for every $x\in E$, then $X\setminus E$ must have positive measure.

Let $\eta_1, \eta_2,\dots$ be a sequence of Borel sections of $\{H_x: x\in X\}$ such 
that $\{\eta_1(x),\eta_2(x),\dots\}$ is dense in the unit ball of $H_x$ for all $x\in X$,   
and consider the sequence of functions $F_k:\mathcal D\to[0,\infty)$, $k\geq 1$,  defined by 
$$
F_k(x,\psi)=|\langle\psi(s)\xi(x),\eta_k(x)\rangle|, \qquad \psi\in \mathcal D_x,\quad x\in X.  
$$
The $F_k$ are Borel functions, they restrict to continuous 
functions on each fiber $\mathcal D_x$, and they have the property that 
$\phi_x$ is maximal on $(s,\xi(x))$ iff 
$$
F_k(x,\psi)\leq \|\phi_x(s)\xi(x)\|, \qquad k=1,2,\dots, \quad \psi\in \mathcal D_x.  
$$  

Consider the following subset of $\mathcal D$
$$
\mathcal D_+=\bigcup_{k=1}^\infty\{(x,\psi)\in \mathcal D: F_k(x,\psi)>\|\phi_x(s)\xi(x)\|\}.  
$$
Being exhibited as 
a countable union of Borel sets, $\mathcal D_+$ is a Borel subset of $\mathcal D$, 
and the natural projection $p: \mathcal D\to X$ restricts to a Borel map 
of $\mathcal D_+$ into $X$.  Let $X_+=p(\mathcal D_+)$ be the range of the restricted 
map.   While $X_+$ is not necessarily a Borel set, Theorem \ref{mtThm3} implies 
that it is an analytic subset 
of $X$ and is therefore absolutely measurable.   Notice too that, 
by its definition,  
$X_+$ is the set of all $x\in X$ such that the UCP map 
$\phi_x:S\to\mathcal B(H_x)$ is {\em not} maximal 
on $(s,\xi(x))$.  

Since $X_+$ is 
absolutely measurable, there exist Borel sets $E_\mu\subseteq X_+\subseteq F_\mu$ 
such that $\mu(F_\mu\setminus E_\mu)=0$.  Note that $\phi_x$ is 
maximal at $(s,\xi(x))$ for every $x\in X\setminus F_\mu\subseteq X\setminus X_+$.  
So by our choice of $s$ and $\xi$, 
$F_\mu$ must have positive measure.   
Hence $E_\mu$ is a Borel subset of $X_+$ of the same positive measure, and we may 
consider the restricted map $p:\mathcal D_0\to E_\mu$ defined by 
$$
\mathcal D_0=\{(x,\psi)\in \mathcal D_+: x\in E_\mu\},\qquad p(x,\psi)=x.  
$$

The projection $p: \mathcal D_0\to E_\mu$ is a surjective Borel function and 
both $\mathcal D_0$ and $E_\mu$ are standard Borel spaces.    
By Theorem \ref{mtThm2} 
there is a Borel set $N\subseteq E_\mu$ of $\mu$-measure zero and a Borel 
cross section $x\in E_\mu\setminus N\mapsto \psi_x\in \mathcal D_x$ for $p$.  
This function $x\mapsto\psi_x$ satisfies (\ref{dbEq2}) for $X_0=E_\mu\setminus N$.  
\end{proof}

\section{Structure of maps with the unique extension property}\label{S:dc}

Given a UCP map $\phi: S\to\mathcal B(H)$ with the unique extension 
property, it is easy to show that 
for every decomposition $\phi=\phi_1\oplus \phi_2$ of $\phi$ into a direct sum 
of UCP maps of $S$, both $\phi_1$ and $\phi_2$ inherit the unique extension 
property.  The purpose of this section is to generalize that fact to infinite 
continuous decompositions of $\phi$ into a direct integral.  
Throughout this section, $(X,\mu)$ will denote a standard Borel 
probability space.  

\begin{thm}\label{dcThm1}
Let 
$\{H_x: x\in X\}$ be a Borel family of separable Hilbert spaces over $X$, 
let $\phi=\{\phi_x: S\to\mathcal B(H_x):x\in X\}$ be a family of UCP maps, 
and let $\mu$ be a probability measure on $X$.  Let 
$$
H=\int_X^\oplus H_x\,d\mu(x)
$$
be the direct integral of Hilbert spaces and let $\phi: S\to \mathcal B(H)$ 
be 
the direct integral of UCP maps 
$$
\phi(s)=\int_X^\oplus \phi_x(s)\,d\mu(x), \qquad s\in S.   
$$
If $\phi$ has the unique extension property, then there is a Borel set $N\subseteq X$ 
of measure zero such that $\phi_x:S\to\mathcal B(H_x)$ has the unique extension property 
for every $x$ in $X\setminus N$.  
\end{thm}

\begin{proof}  Contrapositively, assume that the direct integral $\phi$ 
has the unique extension property but that there is no Borel set 
$N$ with $\mu(N)=0$ such that 
$\phi=\{\phi_x:S\to\mathcal B(H_x)\}$ 
has the unique extension property for every $x\in X\setminus N$.

By Theorem \ref{dbThm1}, there is a Borel set $X_0\subseteq X$ of 
positive measure and a Borel section 
$\psi: x\in X_0\to \psi_x\in UCP(S,\mathcal B(H_x\oplus\ell^2))$ that dilates 
the restricted family 
$\{\phi_x:S\to\mathcal B(H_x): x\in X_0\}$ nontrivially for every $x\in X_0$.   
Thus, for every $x\in X_0$, $\psi_x$ is a dilation of $\phi_x$ that 
cannot be decomposed into a direct sum $\phi_x\oplus\lambda_x$.  
Define a larger bundle of separable Hilbert spaces $\{K_x: x\in X\}$ by 
$$
K_x=
\begin{cases}
H_x\oplus\ell^2,& x\in X_0,\\
H_x,& x\in X\setminus X_0.  
\end{cases}
$$
and a new measurable family $\tilde \phi=\{\tilde\phi_x:S\to\mathcal B(K_x): x\in X\}$ by 
$$
\tilde\phi_x=
\begin{cases}
\psi_x,& x\in X_0,\\
\phi_x,& x\in X\setminus X_0.  
\end{cases}
$$
The family $\{\tilde\phi_x:S\to\mathcal B(K_x): x\in X\}$ is a dilation of the 
original family 
$\{\phi_x:S\to \mathcal B(H_x): x\in X\}$, and we can form the 
direct integral of UCP maps of $S$ 
$$
\tilde\phi=\int_X^\oplus \tilde\phi_x\,d\mu(x)  
$$
acting on the Hilbert space 
$$
K=\int_X^\oplus K_x\,d\mu(x).  
$$
For each $x\in X$, let $P_x$ be the projection of $K_x$ 
on $H_x$, so that the projection $P$ of $K$ on $H$ decomposes into a direct 
integral 
$$
P=\int_X^\oplus P_x\,d\mu(x).  
$$
The UCP map $\tilde\phi$ is a dilation of $\phi$, and $\phi$ is maximal since it 
has the unique extension property.  Hence $P$ commutes with $\tilde\phi(S)$. 
Let $\{a_1,a_2,\dots\}$ be a countable norm-dense subset of $S$.   
Since $P$ is a decomposable operator that
commutes with the sequence of decomposable operators 
$\{\tilde\phi(a_1), \tilde\phi(a_2),\dots\}$, it follows that 
$P_x$ must commute with $\{\tilde\phi_x(a_1), \tilde\phi_x(a_2),\dots\}$ for all $x$ in the 
complement of some $\mu$-null Borel set $N\subseteq X$, and hence $P_x$ commutes 
with $\tilde\phi_x(S)$ for all $x\in X\setminus N$.   
In particular, for $x\in X_0\setminus N$, this implies that 
the constructed family of dilations $\psi_x$ decomposes into a direct sum 
of UCP maps 
$$
\psi_x=\phi_x\oplus\lambda_x, 
$$  
contradicting the stated property of $\psi$ on a set of positive measure.  
\end{proof}

\section{Existence of boundary representations}\label{S:xb}

In this section we prove the following main result.

\begin{thm}\label{xbThm1}
Every separable operator system $S\subseteq C^*(S)$ 
has sufficiently many boundary representations.  
\end{thm}

\begin{proof}  We will show that there is a set $\{\sigma_x: x\in A\}$ of boundary 
representations that satisfies (\ref{xbEq0}).  
To that end, we realize $S\subseteq \mathcal B(H_0)$ as an operator system acting on a separable 
Hilbert space $H_0$.  By Theorem \ref{suThm1}, there is a separable Hilbert space 
$H\supseteq H_0$ and a UCP map $\phi: S\to \mathcal B(H)$ with the 
unique extension property 
such that $P_{H_0}\phi(s)\restriction_{H_0}=s$, $s\in S$.  Obviously, 
$\phi$ is a complete isomorphism of $S$ onto $\phi(S)\subseteq \mathcal B(H)$.  

Choose a maximal abelian von Neumann subalgebra $\mathcal M$ of the 
commutant $\phi(S)^\prime$.  Since $\mathcal M$ acts on a separable Hilbert space, 
it contains a separable unital \cstar\   $\mathcal A$ whose weak closure 
is $\mathcal M$.  The Gelfand spectrum $X$ of $\mathcal A$ is a compact 
metrizable space, hence we may view it 
as a standard Borel space, and there is a probability 
measure $\mu$ on $X$ such that $\mathcal M\cong L^\infty(X,\mu)$.  

Conventional multiplicity theory implies that there is 
a standard Hilbert bundle 
$\{H_x: x\in X\}$ that gives rise to a decomposition
$$
H=\int_X^\oplus H_x\,d\mu(x)
$$
in such a way that $\mathcal M$ is realized as 
the algebra of multiplications 
by scalar functions in $L^\infty(X,\mu)$.  
For example,  this encapsulates the discussion on page 55 of \cite{arvInv}.  
In these  
``coordinates", $\phi$ becomes a UCP map of $S$ into $\mathcal B(H)$ whose 
range $\phi(S)\subseteq \mathcal M^\prime$ consists of decomposable operators.  
The \cstar\ generated by $\phi(S)$ is separable, hence 
the \cstar\ $\mathcal B$ generated by $\mathcal A\cup \phi(S)$ is a 
separable \cstar\ contained in the commutant of $\mathcal M$.  Note that 
by the choice of $\mathcal A$, 
the commutant of $\mathcal B$ is $\mathcal M$, so 
by the double commutant theorem, $\mathcal B$ is weak$^*$-dense in 
the von Neumann algebra $\mathcal M^\prime$ of all decomposable operators. 

Corollary 2 of Theorem 4.2.1 of \cite{arvInv}, together with the 
``flattening and piecing together" procedure described in Section \ref{S:fa},  
implies that there 
is a Borel-measurable family of representations 
$\pi_x: \mathcal B\to \mathcal B(H_x)$, $x\in X$,  with the property 
$$
b=\int^\oplus_X \pi_x(b)\,d\mu(x),\qquad b\in \mathcal B.  
$$ 
Since $\mathcal B$ is weak$^*$- dense in $\mathcal M^\prime$, 
Corollary 2 of Proposition 4.2.2 of \cite{arvInv} implies that $\pi_x(\mathcal B)$ is an 
irreducible set of operators for almost every $x\in X$.  By discarding 
a Borel set of measure zero from $X$, we can assume that $\pi_x(\mathcal B)$ is 
an irreducible \cstar\ for every $x\in X$.  

We can now define a family of UCP maps $\{\phi_x: S\to\mathcal B(H_x): x\in X\}$ by 
setting $\phi_x(a)=\pi_x(\phi(a))$, $a\in S$, $x\in X$, thereby obtaining 
a disintegration of 
$\phi$ into a direct integral 
of UCP maps 
\begin{equation}\label{xbEq1}
\phi(a)=\int_X^\oplus \phi_x(a)\,d\mu(x), \qquad a\in S.   
\end{equation}
After noting that for every $x\in X$, 
$$
\pi_x(\mathcal B)=\pi_x(C^*(\phi(S)\cup\mathcal A))=
C^*(\pi_x(\phi(S)\cup\mathcal A))=C^*(\phi_x(S)),
$$ 
one finds that $\phi_x(S)$ is an irreducible operator system for every $x\in X$.  
Thus, $\{\phi_x: S\to\mathcal B(H_x): x\in X\}$ defines a (Borel) family of irreducible 
UCP maps of $S$.  The decomposition (\ref{xbEq1}) itself implies that 
\begin{equation*}
\|a\|=\|\phi(a)\|=\esup_{x\in X}\|\phi_x(a)\|, \qquad a\in S, 
\end{equation*}
$\esup$ denoting the essential supremum with respect 
to the measure $\mu$.  Note that 
similar formulas hold throughout the matrix hierarchy over $S$.  
Indeed,  for each $n=2,3,\dots$, the direct sum $n\cdot H$ of $n$ copies of $H$ 
decomposes into a direct integral 
$$
n\cdot H=\int_X^\oplus n\cdot H_x\,d\mu(x)
$$ 
and the 
associated map of $n\times n$ matrices 
$(a_{ij})\in M_n(S)\mapsto (\phi(a_{ij}))\in\mathcal B(n\cdot H)$ 
admits a similar direct integral decomposition in 
which the $n\times n$ operator matrix $(\phi(a_{ij}))$ is realized 
as a direct integral of $n\times n$ operator matrices over $\mathcal B(H_x)$, $x\in X$ 
$$
(\phi(a_{ij}))=\int_X ^\oplus (\phi_x(a_{ij}))\,d\mu(x).  
$$
From this formula one concludes that for every $n=2,3,\dots$, 
\begin{equation}\label{xbEq2}
\|(a_{ij})\|=\|(\phi(a_{ij}))\|=\esup_{x\in X}\|(\phi_x(a_{ij}))\|, \quad (a_{ij})\in M_n(S).  
\end{equation}

By Theorem \ref{dcThm1}, there is a Borel set $N\subseteq X$ of measure 
zero such that for each $x\in X\setminus N$, 
$\phi_x$ has the unique extension property, and therefore 
defines an element $\sigma_x$ of $\partial_S$.  Since $S$ is separable, 
we can restrict all terms $a_{ij}$ appearing in the formulas of (\ref{xbEq2}) 
to a countable 
dense subset $C$ of $S$, thereby obtaining  an equivalent countable 
set of formulas.  Finally, since the essential supremum 
of an $L^\infty$ function defined on a measure space is the pointwise supremum 
of that function restricted to a subset whose complement is of measure 
zero, there is a single Borel set $A\subseteq X$, such that $\mu(X\setminus A)=0$, 
with the property 
that the entire sequence of formulas (\ref{xbEq2}) holds as pointwise suprema 
- over $x\in A$ - of norms of  
matrix functions whose entries $\phi_x(a_{ij})$ involve terms with   
$a_{ij}\in C$.  
Finally, since $C$ is dense in $S$, it follows that the set of boundary representations 
$\{\sigma_x: x\in A\}$ satisfies (\ref{xbEq0}).  
\end{proof}

\section{Pure states of $S$}\label{S:ps}

In this section we sharpen Theorem \ref{xbThm1} by showing 
that every pure state of $S$ can be associated with a boundary representation.

\begin{defn}
A state $\rho$ of $C^*(S)$ is called an $S$-{\em boundary} state 
if it is a pure state of $C^*(S)$ and 
the irreducible representation $\pi$ occurring in its GNS representation 
\begin{equation}\label{psEq1}
\rho(x)=\langle\pi(x)\xi,\xi\rangle,\qquad x\in C^*(S) 
\end{equation}
is a boundary representation for $S$.  
\end{defn}

By a  {\em state} of $S$ we mean a positive linear functional $\rho$ on $S$
satisfying $\rho(\mathbf 1)=1$.  
A {\em pure state} of $S$ is  an extreme point of the 
convex weak$^*$-compact set of all 
states of $S$.

\begin{thm}\label{psThm1}
Every pure state of a separable operator system $S$ 
can be extended to an $S$-boundary state of $C^*(S)$.  
\end{thm}

Equivalently, the assertion is that every pure state $\rho$  of $S$ can 
be written the form 
(\ref{psEq1}) $\rho(a)=\langle\pi(a)\xi,\xi\rangle$, $a\in S$, 
for some $\pi\in\partial_S$.  The proof requires the following 
measure-theoretic refinement of the notion of extreme point:

\begin{lem}\label{psLem1}
Let $\phi$ be a pure state of $S$ and let $(X,\mu)$ be a standard Borel 
probability space.  For every  $x\in X$,  let $\rho_x$ be a state  
of $S$ such that for every $a\in S$, $\rho_x(a)$ is a Borel-measurable 
function of $x$, satisfying  
$$
\phi(a)=\int_X\rho_x(a)\,d\mu(x), \qquad a\in S.  
$$
Then $N=\{x\in X: \rho_x\neq\phi\}$ is a Borel set of measure zero.  
\end{lem}

\begin{proof}
Since $S$ is separable, its state space $Y$ is a compact convex metrizable 
space relative to its weak$^*$ topology, 
and $x\in X\mapsto \rho_x\in Y$ is a Borel map of $X$ into $Y$.  Let $\nu$ be 
the push-forward of $\mu$, i.e., the probability measure on $Y$ defined 
on Borel sets by 
$$
\nu(E)=\mu\{x\in X: \rho_x\in E\}, \qquad E\subseteq Y.  
$$
By the standard change-of-variables formula of measure theory, for every bounded Borel 
function $F:Y\to \mathbb C$ we have 
$$
\int_XF(\rho_x)\,d\mu(x)=\int_Y F(\rho)\,d\nu(\rho).  
$$
Taking $F(\rho)=\rho(a)$ for fixed $a\in S$, we obtain
$$
\int_X\rho_x(a)\,d\mu(x)=\int_Y \rho(a)\,d\nu(\rho),   
$$
and hence 
$$
\phi(a)=\int_Y\rho(a)\,d\nu(\rho),\qquad a\in S.  
$$
A result of Bauer (see Proposition 1.4 of \cite{phelps} or \cite{phelps2}) implies that $\nu$ is 
the point mass concentrated at $\phi$.   Let 
$N=\{x\in X: \rho_x\neq \phi\}$.  Then $N$ is a Borel subset 
of $X$ such 
that 
$
\mu(N)=\nu(Y\setminus\{\phi\})=0. 
$    
\end{proof}

\begin{proof}[Proof of Theorem \ref{psThm1}]
Let $\rho$ be an extension of $\phi$ to a state 
of $C^*(S)$, and let $\rho(x)=\langle\pi_0(x)\xi,\xi\rangle$, 
$x\in C^*(S)$, be its GNS representation.  Note that the Hilbert 
space $[\pi_0(C^*(S))\xi]$ of $\pi_0$ is separable.  
By Theorem \ref{suThm1}, the UCP map $\pi_0\restriction_S$
can be dilated to a UCP map of $S$, on a larger 
separable Hilbert space, which has the unique extension property.  
Let $\pi$ be the extension of this dilation to $C^*(S)$.  
The formula $\phi(a)=\langle\pi(a)\xi,\xi\rangle$ persists 
for $a\in S$. 

As in the proof of Theorem \ref{xbThm1}, we can 
decompose $\pi$ into a direct integral of irreducible representations 
$\pi_x:C^*(S)\to\mathcal B(H_x)$ parameterized by a standard Borel probability space $(X,\mu)$.  
By Theorem 
\ref{dcThm1}, $\pi_x$ is a boundary representation for almost every $x\in X$.  The vector 
$\xi$ becomes a square-integrable section $x\in X\mapsto \xi(x)\in H_x$, 
and we can define 
a Borel section of unit vectors over the Borel set $X_0=\{x\in X: \xi(x)\neq 0\}$ by
$$
e(x)=
\|\xi(x)\|^{-1}\xi(x),\qquad x\in X_0.  
$$
This decomposition of $\pi$ leads to the following representation of $\phi$
\begin{align*}
\phi(a)&=\int_X\langle\pi_x(a)\xi(x),\xi(x)\rangle\,d\mu(x)
\\
&=
\int_X\langle\pi_x(a)e(x),e(x)\rangle\|\xi(x)\|^2\,d\mu(x)
=\int_{X_0}\langle\pi_x(a)e(x),e(x)\rangle\,d\nu(x),  
\end{align*}
for $a\in S$, where $\nu$ is the probability measure defined 
on $X_0$ by 
$$
d\nu(x)=\|\xi(x)\|^2\,d\mu(x).
$$  
Note that $\rho_x(b)=\langle\pi_x(b)e(x),e(x)\rangle$ is an $S$-boundary state 
of $C^*(S)$ for every $x\in X_0$.  
Since $\phi$ is a pure state of $S$, 
 Lemma \ref{psLem1} implies that 
$$
\phi(a)=\rho_x(a)=\langle\pi_x(a)e(x),e(x)\rangle, \qquad a\in S, 
$$
for all $x$ in the complement $X_0\setminus N$ of a Borel set of $\nu$-measure zero,  
thereby exhibiting many  $S$-boundary states that extend $\phi$.  
\end{proof}

\section{Concluding remarks}\label{S:cr}

 The proof of Theorem \ref{xbThm1} is far from constructive.  Rather, 
it is more akin to probabilistic arguments whereby one 
establishes the existence of 
a desired property by constructing a nonvacuous probability space 
in which the property can be shown to hold almost surely.  The 
proof of Theorem \ref{psThm1} illustrates the technique in this context.  

Naturally, it would be desirable to 
get rid of the disintegration theory that is seriously exploited 
above by finding a more direct construction of 
boundary representations.  A preliminary attempt to do that was made in 
\cite{arvSubalgI}, but without much success.  
Indeed, it is still unclear how one might 
effectively characterize the pure states of $C^*(S)$ whose 
GNS representations are boundary representations for $S$.  For example, 
can every pure state $\rho$ of $S$ be extended to a pure state $\tilde\rho$ of 
$C^*(S)$ whose GNS representation is a boundary representation for $S$?  
The answer is yes if $S$ is separable by Theorem \ref{psThm1}, 
or if $C^*(S)$ is commutative in general.    
What we are proposing is a more direct proof that will work for 
inseparable operator systems.  

As a test problem for such developments, we propose:

\vskip0.1in
\noindent
{\bf Problem:} Does Theorem \ref{xbThm1} remain true for inseparable operator systems?  
\vskip0.1in

Perhaps it is worth pointing out that in general, heroic 
attempts to get rid of separability hypotheses for problems 
in operator algebras can force 
one to look carefully at the fundamentals of set theory.  
For example,  Akemann and Weaver  \cite{akW}
have constructed a counter example to Naimark's problem 
by making use of a set-theoretic principle that is 
known to be consistent with,  but not provable from,  the standard axioms 
of set theory.  They also showed that the statement 
{\em There is a counter example to Naimark's problem that is generated 
by $\aleph_1$ elements} is undecidable within standard set theory.  


\bibliographystyle{alpha}

\newcommand{\noopsort}[1]{} \newcommand{\printfirst}[2]{#1}
  \newcommand{\singleletter}[1]{#1} \newcommand{\switchargs}[2]{#2#1}

\end{document}